\documentclass[12pt,reqno]{amsart}
  \usepackage[all,2cell,arrow]{xy}
  \usepackage{amsfonts}
  \usepackage{amsthm}
  \usepackage{amsmath}
  \usepackage{amssymb}
\usepackage{mathtools}
\usepackage{tensor}
 \numberwithin{equation}{section}
 \usepackage{xcolor}
\usepackage{graphicx}
\usepackage[pagebackref,colorlinks,citecolor=darkgreen,linkcolor=darkgreen]{hyperref}

\usepackage{enumerate,xspace}

\UseAllTwocells
\CompileMatrices

 \def\black{\color{black}}

 \definecolor{darkgreen}{rgb}{0,0.45,0} 
\renewcommand{\epsilon}{\varepsilon}
\renewcommand{\phi}{\varphi}

\newcommand{\ca}{\ensuremath{\mathcal A}\xspace}

\newcommand{\ce}{\ensuremath{\mathcal E}\xspace}

\newcommand{\cj}{\ensuremath{\mathcal J}\xspace}
\newcommand{\ck}{\ensuremath{\mathcal K}\xspace}
\newcommand{\cl}{\ensuremath{\mathcal L}\xspace}
\newcommand{\cm}{\ensuremath{\mathcal M}\xspace}

\newcommand{\cp}{\ensuremath{\mathcal P}\xspace}

\newcommand{\cv}{\ensuremath{\mathcal V}\xspace}

\newcommand{\bbi}{\ensuremath{\mathbb I}\xspace}

\newcommand{\atwo}{{\mathbf 2}}

\newcommand{\id}{\textnormal{id}}

\newcommand{\Enr}{\ensuremath{\textnormal{-}\Cat}\xspace}

\newcommand{\Cat}{\ensuremath{\mathbf{Cat}}\xspace}

\newcommand{\Set}{\ensuremath{\mathbf{Set}}\xspace}
\newcommand{\SSet}{\ensuremath{\mathbf{SSet}}\xspace}

\newcommand{\LP}{\textnormal{LP}}

\newcommand{\iso}{\bbi}

\newcommand{\op}{\ensuremath{^{\textnormal{op}}}}

\DeclareMathOperator{\cod}{cod}
\DeclareMathOperator{\dom}{dom}

\DeclareMathOperator{\TF}{\mathbf{TF}}

\DeclareMathOperator{\colim}{colim}

\DeclareMathOperator{\Ps}{Ps}

\newcommand{\two}{\ensuremath{\mathbf{2}\xspace}}

\newcommand{\x}{\times}

\def\1c#1{\stackrel{#1}{\to}}

  \newtheorem{proposition}{Proposition}[section]
  \newtheorem{lemma}[proposition]{Lemma}
  
  \newtheorem{corollary}[proposition]{Corollary}
  
  \newtheorem{theorem}[proposition]{Theorem}
  \newtheorem{defprop}[proposition]{Definition/Proposition}

  \theoremstyle{definition}
  \newtheorem{definition}[proposition]{Definition}

  \newtheorem{example}[proposition]{Example}
  \newtheorem{examples}[proposition]{Examples}

  \theoremstyle{remark}
  \newtheorem{remark}[proposition]{Remark}

  \newcounter{c}
  \renewcommand{\[}{\setcounter{c}{1}$$}
  \newcommand{\etyk}[1]{\vspace{-7.4mm}$$\begin{equation}\Label{#1}
  \addtocounter{c}{1}}
  \renewcommand{\]}{\ifnum \value{c}=1 $$\else \end{equation}\fi}
\newcommand\isf{{\mathrel{\rotatebox{270}{\scriptsize$\twoheadrightarrow$}}}}

\begin{document}

\title[On $2$-categorical $\infty$-cosmoi]{On $2$-categorical $\infty$-cosmoi}
\author{John Bourke}
\address{Department of Mathematics and Statistics, Masaryk University, Kotl\'a\v rsk\'a 2, Brno 61137, Czech Republic}
\email{bourkej@math.muni.cz}

\author{Stephen Lack}
\address{School of Mathematical and Physical Sciences, Macquarie
 University NSW 2109, Australia}
\email{steve.lack@mq.edu.au}
\thanks{The first-named author acknowledges the support of the Grant Agency of the Czech Republic under the grant  22-02964S. The second-named author acknowledges with gratitude the support of an
 Australian Research Council Discovery Project DP190102432.}

\date{\today}
\maketitle

\begin{abstract}
Recently Riehl and Verity have introduced $\infty$-cosmoi, which are
certain simplicially enriched categories with additional structure.
In this paper we investigate those $\infty$-cosmoi which are  in fact $2$-categories; we shall refer to these as $2$-cosmoi.  

We show that each $2$-category with flexible limits gives rise to a $2$-cosmos whose distinguished class of isofibrations consists of the normal isofibrations.  Many examples arise in this way, and we show that such $2$-cosmoi are minimal as Cauchy-complete $2$-cosmoi.  Finally, we investigate accessible $2$-cosmoi and develop a few aspects of their basic theory.
\end{abstract}

\section{Introduction}

In recent years Riehl and Verity have introduced the notion of an
$\infty$-cosmos as a framework in which to do $(\infty,1)$-category
theory in a model independent way \cite{elements}.  An $\infty$-cosmos
is a simplicially-enriched category $\ck$ equipped with a class of
morphisms, called isofibrations, satisfying certain properties,
notably the existence of certain basic kinds of limits.  Using these
limits, one can develop much of the theory of $\infty$-categories
inside an $\infty$-cosmos, such as the study of $\infty$-categories
with limits and/or colimits of a given shape, or cartesian fibrations of $\infty$-categories, and these results then apply to various models of $\infty$-categories, such as quasicategories or complete Segal spaces.

The study of $\infty$-categories using $\infty$-cosmoi is very much $2$-categorical in nature.  In particular, each $\infty$-cosmos $\ck$ has an associated \emph{homotopy $2$-category} $h\ck$ and many of the definitions in an $\infty$-cosmos, such as equivalences and adjunctions, take place within this homotopy $2$-category \cite{Riehl2015The-2-category}.  This means that $\infty$-categorical results such as ``right adjoints preserve limits" can be given $2$-categorical proofs.  At the same time, much of the theory of $\infty$-cosmoi is guided by earlier results in $2$-category theory such as Street's study of fibrations in a $2$-category \cite{Street1972Fibrations} and the theory of flexible limits \cite{BKPS}.

In the opposite direction, $2$-categories can be viewed as
simplicially-enriched categories, by taking the nerves of their
hom-categories, 
 Therefore it makes sense to ask whether naturally
occuring $2$-categories, such as the 2-category of monoidal categories
and strong monoidal functors, or that of coherent categories and coherent
functors, underlie $\infty$-cosmoi.  The main goal of this paper is to
give a positive general answer to this question.  In particular, in
Theorem~\ref{thm:flexible} we prove that each $2$-category with
flexible limits (in the sense of Bird et al \cite{BKPS}) forms part of
an $\infty$-cosmos, when equipped with the class of \emph{normal
  isofibrations}.  As explained in \cite{BKPS}, many $2$-categories of
categorical structures and their \emph{pseudomorphisms} have flexible
limits, so that this is a rather general result.

We use the term \emph{$2$-cosmoi} to refer to these $2$-categorical examples of $\infty$-cosmoi --- namely, $2$-categories equipped with the structure of an $\infty$-cosmos on the associated simplicially-enriched category.
For an $\infty$-cosmos $\ck$ which arises in this way from a 2-category,  the homotopy 2-category $h\ck$ is just the
original 2-category. Thus any $\infty$-cosmological notions which are
defined in the homotopy 2-category, such as those discussed in \cite[Chapter~2]{elements}, agree with the usual 2-categorical
notions in this case. 

We say that a 2-cosmos is {\em Cauchy complete} if it is Cauchy
complete as a 2-category, meaning that idempotents split, and if moreover  the isofibrations are closed under retracts, as is the
case in the key example of a flexibly complete 2-category equipped with the normal
isofibrations. And in fact this notion of Cauchy complete 2-cosmos
shows that the choice of the normal isofibrations is not just widely
applicable but natural: we show in Theorem~\ref{thm:minimal}
that for any Cauchy complete 2-cosmos, the chosen
class of isofibrations must contain the normal isofibrations. Thus the
normal isofibrations are the {\em minimal} choice, at least in the Cauchy
complete context.

As well as Cauchy completeness for 2-cosmoi, we also study
accessibility: this is the $2$-categorical version of the
accessible $\infty$-cosmoi introduced in \cite{Bourke2021Accessible}. 

Let us now give an overview of the paper. In
Section~\ref{sect:basic}, we review the background and give an
elementary characterization of $2$-cosmoi, involving only
$2$-categorical rather than simplicially-enriched concepts.  In
Section~\ref{sect:isofibrations}, we investigate the different
flavours of isofibration that naturally arise in $2$-category theory,
with a particular focus on the normal isofibrations.  In
Section~\ref{sect:flex} we describe the $\infty$-cosmos structure on a
$2$-category with flexible limits and provide a host of examples.
Finally, in Section~\ref{sect:accessible}, we turn to accessible
$2$-cosmoi.

\section{Basic facts about $\infty$-cosmoi and $2$-cosmoi}
\label{sect:basic}

Let $\Cat$ denote the cartesian closed category of small categories, and let $\SSet$ denote the cartesian closed category of simplicial sets.  Central to our concerns is the adjunction

\begin{equation*} 
\xymatrix{
\Cat \ar@{}[r]|-{\bot} \ar@<-4.5pt>[r]_-{{N}} & \SSet
  \ar@<-4.5pt>[l]_-{\Pi}}
\end{equation*}
whose right adjoint $N$ takes the nerve of a small category and whose
left adjoint $\Pi$ takes the classifying category of a simplicial set.
It is well known that the left adjoint $\Pi$, as well as the right
adjoint $N$, preserves finite products.  It follows that the
adjunction lifts to a further adjunction

\begin{equation*} 
\xymatrix{
\Cat\Enr \ar@{}[r]|-{\bot} \ar@<-4.5pt>[r]_-{{N_*}} & \SSet\Enr
  \ar@<-4.5pt>[l]_-{\Pi_*}}
\end{equation*}
between the categories of $\Cat$-enriched categories
(a.k.a. $2$-categories) and of simplicially-enriched categories.  When
$\ck$ is an $\infty$-cosmos, $\Pi_*\ck$ is the homotopy $2$-category
$h\ck$ mentioned in the introduction, and which we will have no need
to discuss further. 

For $\ck$ a $2$-category, $N_*\ck$ has the same objects as $\ck$ and
homs $N_*\ck(A,B) = N(\ck(A,B))$ obtained by taking the nerve.  
 We write $n\ck$ as an abbreviation for $N_*\ck$.  Our goal is now to investigate what it means to equip $n\ck$ with the structure of an $\infty$-cosmos.  

First, let us recall from \cite{elements} that an {\em $\infty$-cosmos} is a simplicially enriched category $\ck$
together with a class of morphisms called \emph{isofibrations}, denoted $A \twoheadrightarrow B$, closed
under composition and containing the isomorphisms, such that 
\begin{enumerate}
\item each hom $\ck(A,B)$ is a quasicategory;
\item the morphism $\ck(A,p)\colon\ck(A,B)\twoheadrightarrow\ck(A,C)$
  is an isofibration of quasicategories for each isofibration $p
  \colon B \twoheadrightarrow C$;
\item \ck has products, powers by simplicial sets, pullbacks along
  isofibrations, and limits of countable towers of isofibrations;
\item the class of isofibrations is stable under pullback; closed under products, limits of countable towers, and
  Leibniz powers by monomorphisms of simplicial sets; and contains all
  maps with terminal codomain.  
\end{enumerate}

\begin{definition}
A \emph{$2$-cosmos} is a $2$-category $\ck$ together with an $\infty$-cosmos structure on $n\ck$.
\end{definition}

As observed in the introduction, the homotopy 2-category of
$n\ck$ is then just $\ck$ itself.

The following proposition describes $2$-cosmoi in purely
$2$-categorical terms.  In the statement, recall that a functor
$f\colon A \to B$ is said to be an \emph{isofibration of categories} if it has the
isomorphism lifting property: namely, given an isomorphism $\alpha
\colon X \to fY$ in $B$, there exists an isomorphism $\alpha'\colon X' \to Y$ in $A$ such that $f\alpha'=\alpha$.

\begin{proposition}\label{prop:2cosmos}
A $2$-cosmos is a $2$-category $\ck$ together with a class of morphisms called \emph{isofibrations}, denoted $A \twoheadrightarrow B$, closed
under composition and containing the isomorphisms, such that 
\begin{enumerate}
\item[(a)] the morphism $\ck(A,p)\colon\ck(A,B)\twoheadrightarrow\ck(A,C)$
  is an isofibration of categories for each isofibration $p
  \colon B \twoheadrightarrow C$;
\item[(b)] \ck has products, powers by small categories, pullbacks along
  isofibrations, and limits of countable towers of isofibrations;
\item[(c)] the class of isofibrations is stable under pullback; closed under products, limits of countable towers, and
  Leibniz powers by injective on objects functors; and contains all
  maps with terminal codomain.  
\end{enumerate}
\end{proposition}

Before proving this result, we will require a couple of preliminary results.  
The first is very well-known:

\begin{lemma}\label{lemma:powers}
There is an isomorphism $[X,NY] \cong N[\Pi X,Y]$, natural in simplicial sets $X$ and categories $Y$.
\end{lemma}
\begin{proof}
By Yoneda, this is equivalent to the statement that we have a bijection $$\SSet(Z,[X,NY]) \cong \SSet(Z,N[\Pi X,Y])$$
natural in all variables.  Now we have a natural bijection $$\SSet(Z,[X,NY]) \cong \SSet(Z \times X,NY) \cong \Cat(\Pi(Z \times X),Y)$$ and similarly $$\SSet(Z,N[\Pi X,Y]) \cong \Cat(\Pi Z, [\Pi X,Y]) \cong \Cat(\Pi Z \times \Pi X,Y)$$ so that we must equally give a natural bijection $$\Cat(\Pi(Z \times X),Y) \cong \Cat(\Pi Z \times \Pi X,Y)$$ but this follows immediately from the fact that $\Pi$ preserves finite products.
\end{proof}

\begin{lemma}\label{lemma:io}
A functor $f\colon A \to B$ is injective on objects if and only if it
is isomorphic in $\Cat^{\atwo}$ to $\Pi(g)$ for some monomorphism $g$ of simplicial sets.
\end{lemma}
\begin{proof}
Since the objects of $\Pi X$ are the $0$-simplices of $X$, the functor
$\Pi$ sends monomorphisms to injective-on-objects functors, giving one
direction.

For the converse, suppose that $f\colon A\to B$ is
injective on objects. We may factorize $Nf \colon NA\to NB$ as 
\[ \xymatrix{
    NA \ar[r]^-i & X \ar[r]^-p & NB } \]
where $i$ is a monomorphism and $p$ is a trivial fibration of
simplicial sets.
Furthermore, since $Nf$ is injective on
$0$-simplices, we may do this in such a way that $p$ is bijective on
$0$-simplices. Applying $\Pi$ gives a factorization $f\cong\Pi(p)\Pi(i)$,
so it will suffice to show that $\Pi(p)$ is invertible.

 Since $NB$ is (the nerve of) a category and $p$ is a trivial fibration, $X$ is a
 quasicategory and $\Pi X$ is the homotopy category of $X$, as
 described in \cite[Lemma~1.1.12]{elements}.
 Since $p$ is bijective on 0-simplices and full on 1-simplices,
 $\Pi(p)$ is bijective on objects and full; thus it will suffice to show
 that $\Pi(p)$ is faithful. 
 Consider then a parallel pair of morphisms in $\Pi(X)$, which we may represent as homotopy classes $[\alpha]$ and $[\beta]$, where
 $\alpha$ and $\beta$ are parallel $1$-simplices in $X$, and suppose that $[p \alpha] = [p \beta]$.  Then $p\alpha$ and $p\beta$ are homotopic in $Y$, and since $p$ is a trivial fibration we obtain
 a lifted homotopy in $X$ between $\alpha$ and $\beta$, and so
 $[\alpha]=[\beta]$, giving faithfulness of $\Pi(p)$. 
\end{proof}

\begin{proof}[Proof of Proposition~\ref{prop:2cosmos}]
We must investigate what it means to equip $n\ck$ with the structure
of an $\infty$-cosmos. To begin with, this involves giving a class of
morphisms in $n\ck$, closed under composition and containing the
isomorphisms.  But since $\ck$ and $n\ck$ have the same underlying
category, this is just a class of morphisms in $\ck$ satisfying the
same two properties.

Condition (i) in the definition of $\infty$-cosmos asserts that $n\ck$ is
enriched in quasicategories.  This is redundant since $\ck$ is
$\Cat$-enriched and the nerve of a category is always a quasicategory.

Condition (ii) asserts that  $N(\ck(A,p))$  is an isofibration of quasicategories for each isofibration $p  \colon B \twoheadrightarrow C$.  By Observation 1.1.19(iv) of \cite{elements} a functor is an isofibration of categories just when its nerve is an isofibration of quasicategories.  Thus Condition~(ii) becomes Condition~(a) above.
  
Next we show that Condition (iii) for $n\ck$ is equivalent to
Condition~(b) for $\ck$.  Condition (iii) says that three kinds of
conical limits, plus powers, exist.  With regards the conical limits,
these are limits in the underlying category $\ck_0$ which have,
moreover, the stronger universal property of being $\SSet$-enriched
limits in $n\ck$ --- that is, limits in $\ck_0$ preserved by each
$n\ck(A,-)\colon \ck_0 \to \SSet$.  Now since $n\ck(A,-)$ is the composite
$N \circ \ck(A,-)\colon \ck_0 \to \Cat \to \SSet$, and since $N$ preserves
and reflects limits, $n\ck(A,-)$ preserves precisely those limits preserved by $\ck(A,-)\colon \ck_0 \to \Cat$ --- that is, limits in $\ck_0$ that are $\Cat$-enriched limits.  

Turning to powers, let $X$ be a simplicial set and $A \in \ck_0$.  We write $A^X_{n\ck}$ for the power of $A$ by $X$ in $n\ck$ and use the corresponding notation for powers in $\ck$.  We claim that we have an isomorphism $A^X_{n\ck} \cong A^{\Pi X}_{\ck}$, either side existing if the other does.    Certainly if $A^X_{n\ck}$ exists, then we have the defining isomorphism $N(\ck(B,A^X_{n\ck})) \cong [X,N(\ck(B,A)]$ natural in $B$.  Composing this with the natural isomorphism $[X,N(\ck(B,A)] \cong N[\Pi X,\ck(B,A)]$ of Lemma~\ref{lemma:powers} and using fully faithfulness of the Yoneda embedding we obtain a natural isomorphism $\ck(B,A^X_{n\ck}) \cong [\Pi X,\ck(B,A)]$, which induces --- now by definition of powers in $\ck$ --- an isomorphism $\phi_{A,X}\colon A^X_{n\ck} \cong A^{\Pi X}_{\ck}$.  The other direction is identical.  Hence $n\ck$ admits powers just when $\ck$ admits powers by small categories of the form $\Pi X$.  But since each small category is of this form, we conclude that $n\ck$ admits powers just when $\ck$ does.  

Finally, we show that Condition (iv) for $n\ck$ becomes Condition (c) for $\ck$.  This is straightforward, bar the Leibniz conditions.  For these, we note that the isomorphisms $\phi_{A,X}\colon A^X_{n\ck} \cong A^{\Pi X}_{\ck}$ constructed above are natural in $A$ and $X$.  Given a monomorphism $i\colon X \to Y$ of simplicial sets and an isofibration $f\colon A \twoheadrightarrow B$, it therefore follows that we have a commutative square

\begin{equation*}
\xymatrix{
A^X_{n\ck} \ar[d]_{\phi_{A,X}}  \ar[r]^{ {i \widetilde\pitchfork f}} & f^i_{n\ck} \ar[d]^{\phi_{i,f}} \\
A^{\Pi X}_{\ck}  \ar[r]^{{(\Pi i) \widetilde\pitchfork f}} & f^{\Pi{i}}_{\ck}
}
\end{equation*}
in which the vertical maps are isomorphisms and the horizontal maps
are the induced ones to the pullback-product, with the upper
horizontal living in $n\ck$ and the lower one in $\ck$.  In
particular, the upper horizontal is an isofibration just when the
lower one is one.  Since, by Lemma~\ref{lemma:io}, the functors of the form $\Pi(i)$ for $i$ mono are, up to isomorphism, the injective on objects functors, we conclude that isofibrations are closed under Leibniz powers by monomorphisms in $n\ck$ just when they are closed under Leibniz powers by injective on objects functors in $\ck$.
\end{proof}

\section{Flavours of isofibration in a $2$-category}
\label{sect:isofibrations}

In this section, we recall the various kinds of isofibrations in a
$2$-category and the relationships between them.

We suppose throughout that $\ck$ is a 2-category with pseudolimits of
arrows, as in Definition/Proposition~\ref{defprop:pseudolimits} below,
and use the notation introduced there.

\subsection{2-categorical preliminaries}

First we recall a few facts and definitions about equivalences in
2-categories:

\begin{defprop}~ \label{defprop:equivalences}
\begin{enumerate}
\item A morphism $f\colon A\to B$ in a 2-category $\ck$ is an {\em
    equivalence} when there exist isomorphisms $\eta\colon 1 \cong
  gf$ and $\epsilon \colon fg\cong 1$; it is then possible to choose
  $\eta$ and $\epsilon$ so that they satisfy the triangle equations.
\item If $f$ is an equivalence and has a section, then it is possible
  to choose the remaining structure so that $\epsilon$ is an identity;
  we then say that $f$ is a {\em retract equivalence}.
\item If $f$ is an equivalence and has a retraction, then it is
  possible to choose the remaining structure so that $\eta$ is an
  identity; we then say that $f$ is an {\em injective equivalence}. 
\end{enumerate}
\end{defprop}

Next we recall a few basic facts about pseudolimits of arrows.

\begin{defprop}~ \label{defprop:pseudolimits}
  \begin{enumerate}
  \item The pseudolimit $L_f$ of an arrow $f\colon A\to B$ is the
    universal diagram as on the left below. Thus for every diagram as
    in the centre, there is a unique arrow $c\colon C\to L_f$ with
    $u_f c=a$, $v_f c=b$, and $\lambda_f c=\beta$. There is also a
    2-dimensional aspect to the universal property.
  \item In particular, there is a unique arrow $d_f\colon A\to L_f$
    with $u_f d_f=1$, $v_f d_f=f$, and $\lambda_f d_f=\id_f$; there is
    also a unique invertible 2-cell $\delta_f\colon 1\cong d_fu_f$
    with $u_f\delta_f=\id$ and $v_f\delta_f=\lambda_f$.
  \item Thus in fact $u_f$ is a retract equivalence and $d_f$ an injective
    equivalence; conversely, if $u_f$ is known to be an equivalence then
    the 2-dimensional aspect of the universal property is automatic.
  \item The pseudolimit $L_{1_A}$ of an identity arrow $1_A\colon A\to
    A$ is equally the power $A^\bbi$ of $A$ by the generic category
    $\bbi$ containing an isomorphism.
  \item The pseudolimit $L_f$ can equally be constructed as a pullback
    as in the square on the right below. There is then a unique
    induced $w_f\colon A^\bbi\to L_f$ making the two triangles
    commute. 
  \end{enumerate}
  
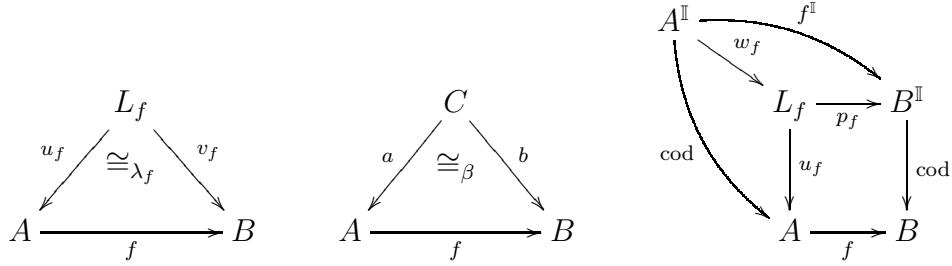
\[ \xymatrix @R1pc {
      &&&&&& A^\bbi \ar@/_1pc/[dddr]_{\cod} \ar@/^1pc/[drr]^{f^\bbi}
      \ar[dr]^{w_f} \\ 
      & L_f \ar[ddl]_{u_f} \ar[ddr]^{v_f}
      \ar@{}[dd]|{\textstyle\cong_{\lambda_f}} &&&
      C \ar[ddl]_a \ar[ddr]^b \ar@{}[dd]|{\textstyle\cong_{\beta}} &&&
      L_f \ar[r]_-{p_f}  \ar[dd]^{u_f} & B^\bbi \ar[dd]^{\cod} \\
      \\ 
      A \ar[rr]_-f && B & A \ar[rr]_-{f} && B && A \ar[r]_-f & B } \]
\end{defprop}

\subsection{Normal isofibrations}

The standard kind of isofibration in a 2-category $\ck$ consists of a morphism $f\colon A\to B$ for which the induced $\ck(C,f)\colon\ck(C,A)\to\ck(C,B)$ is an isofibration of categories for all $C\in \ck$.  In elementary terms, this means that given $g\colon X \to A$, $h\colon X \to B$ and an isomorphism $\alpha\colon f \circ g \cong h$ there exists an isomorphism $\alpha^{\prime}\colon g \cong h^{\prime}$ such that $f \circ \alpha^{\prime} = \alpha$.   We shall call these {\em representable isofibrations}, to distinguish them from the isofibrations that are specified as part of the structure of $\infty$-cosmos.

A \emph{discrete isofibration} is a representable isofibration such that the lifting $(\alpha^\prime, h^\prime)$ above is unique.

A {\em cleavage} for a representable isofibration consists of a
choice, for each $(g,\alpha,h)$, of a pair
$(\alpha^{\prime},h^{\prime})$ as above such that these choices are
\emph{natural in $X$}.  The cleavage is said to be \emph{normal} if
whenever $\alpha$ is an identity $2$-cell so is $\alpha^{\prime}$.  A
\emph{normal isofibration} \cite{Garner:2DTypeTheory} is, by
definition, a representable isofibration supporting a normal cleavage.  

Each discrete isofibration is a normal isofibration.  Moreover in
$\Cat$ each representable isofibration is normal and the corresponding
statement is true in many $2$-categories of categories with structure:
 see also Section~\ref{sect:NIP} below. 

\begin{proposition}\label{prop:characterizations-via-L}
 Let $f\colon A\to B$ be an arrow in $\ck$, and form the pseudolimit
 $L_f$ of $f$. Then
 \begin{enumerate}
 \item $f$ is a representable isofibration if and only if there is an
   invertible 2-cell $\kappa\colon x\cong u_f$ with
   $f\kappa=\lambda_f$ (and so also $fx=v_f$)
 \item $f$ is a normal isofibration if and only if moreover $\kappa\colon
   x\cong u_f$ can be chosen with $\kappa d_f$ an identity.
 \end{enumerate}
\end{proposition}

\begin{proof}
The existence of a lifting $\kappa\colon x\cong u_f$ of $\lambda_f$ as
in (i) is part of the definition of representable
isofibration. Conversely, given $\kappa\colon x\cong u_f$ and a
lifting problem $\beta\colon b\cong fa$, by the universal property of
$L_f$ there is an induced $c$ with $u_fc=a$, $v_fc=b$, and $\lambda_f
c=\beta$ and now $\kappa c\colon xc\cong u_fc=a$ gives the required
cleavage.

Similarly, the existence of $\kappa\colon x\cong u_f$ satisfying the
extra condition in (ii) is part of the definition of normal
isofibration. For the converse, if $\beta\colon b\cong fa$ is an
identity, then the induced $c$ factorizes through $d_f$, and so
$\kappa c$ is also an identity. 
\end{proof}

In Proposition~\ref{prop:wfs} below we describe a weak factorization system on a suitable
2-category $\ck$ whose
right class consists of the normal isofibrations; this implies
many of the stability properties of normal isofibrations that we will
need. It is the (cofibration, trivial fibration) weak factorization
system for the natural model structure \cite[Section~4]{hty2mnd} on $\ck\op$. The
connection with normal isofibrations, as defined above, was made in \cite[Remark~3.4.8]{Garner:2DTypeTheory}, and it also
featured in \cite{Gambino:hty2limits}.

\begin{proposition}\label{prop:wfs}
Let $\ck$ be a $2$-category with pseudolimits of arrows.  Then there is a weak factorization system on $\ck_0$ whose left class consists of the injective equivalences and whose right class consists of the normal isofibrations.  
\end{proposition}
\begin{proof}
The factorization of an arrow $f$ will be given by $f=v_f d_f$ as in
Definition/Proposition~\ref{defprop:equivalences}. For the fact that $d_f$
is an injective equivalence, and that this gives a weak factorization
system consisting of the injective equivalences and {\em some} class $\cp$
of arrows, see \cite[Section~4]{hty2mnd}. In that section, a model
structure is described for any 2-category with finite limits and
finite colimits. If one applies this to $\ck\op$ and considers the
(cofibration, trivial fibration) weak factorization system, one
obtains a weak factorization system (injective equivalence, $\cp$) on
$\ck$. Note that the finite limits and finite colimits are assumed
only because they are part of the definition of a model structure;
only pseudolimits of arrows are required for this weak factorization
system. The dual case of \cite[Section~4]{hty2mnd}, as used here, was
explicity considered in \cite{Gambino:hty2limits}.

We now turn to the fact that the right class $\cp$ consists of the
normal isofibrations. This was observed without proof in
\cite[Remark~3.4.8]{Garner:2DTypeTheory}, so we now give a proof.
First we show that the normal isofibrations have the right lifting
property with respect to the injective equivalences; in other words,
every normal isofibration is in $\cp$. To see this, consider a
commutative square 
\begin{equation*}
\xymatrix{
A \ar[d]_{i} \ar[r]^{f}  & C \ar[d]^{p} \\
B \ar[r]_{g} & D}
\end{equation*}
in which $i$ is an injective equivalence and $p$ is a normal
isofibration. Then there exist $r\colon B \to A$ with $ri=1_A$ and an
invertible $2$-cell $\eta\colon 1_B \cong ir$ satisfying the triangle
equations $\eta i = 1$ and $r\eta = 1$.  

The invertible 2-cell $g\eta\colon g \cong gir = pfr$ lifts along $p$
to give a $1$-cell $h\colon B \to C$ and invertible $2$-cell
$\theta\colon h \cong fr$ satisfying  $ph = g$ and $p\theta = g \eta$.
To prove that $h$ is a diagonal filler, we must show that $hi=f$.  To
this end, observe that by naturality of the cleavage, $hi$ is the
domain of the lifting of $g\eta i\colon gi \cong giri = pfri = pf$ but
now $g\eta i =\id$ and hence, by normality, $hi = f$.

For the converse, suppose that $f\colon A\to B$ is in $\cp$, and form
the factorization $f=v_f d_f$. Then there exists a filler $x$ as in
\[ \xymatrix{
    A \ar[r]^-1 \ar[d]_{d_f} & A \ar[d]^f \\
    L_f \ar[r]_-{v_f} \ar[ur]^x & B. } \]
Now $u_f d_f=1$, and as observed in
Definition/Proposition~\ref{defprop:pseudolimits} as well as in
\cite[Section~4]{hty2mnd},  there is an
isomorphism $\delta_f\colon 1\cong d_f u_f$ satisfying the triangle
equations and $f\delta_f=\lambda_f$. Now $x\delta_f\colon x\cong xd_fu_f=u_f$
satisfies $fx\delta_f=v_f\delta_f=\lambda_f$, and, by one of the
triangle equations $x\delta_f d_f=\id$. Thus $f$ is a normal
isofibration by Proposition~\ref{prop:characterizations-via-L}.
\end{proof}

The following lemma is useful in recognising $2$-cosmoi.

\begin{lemma}\label{lemma:Leibniz}
Let $\ck$ be a $2$-category with products and powers. Let $\cp$ be a
class of morphisms in $\ck$ such that:
\begin{itemize}
\item $\cp$ is closed under products and composition
\item pullbacks of maps in $\cp$ along arbitrary maps exist and are in $\cp$.
\end{itemize}
If $\cp$ contains the discrete isofibrations, then it is closed under 
Leibniz powers by the injective on objects functors. 
\end{lemma}
\begin{proof}
Let $\cj$ be the class of all functors $j$ with the property that
$\cp$ is closed under Leibniz powers by $j$. We need to show that
$\cj$ contains the injective on objects functors.

A Leibniz power by a functor $0\to W$ with initial domain is just an
ordinary power by $W$; thus $\cj$ contains such functors for $W$ a
discrete category. An
arbitrary injective on objects functor $j\colon X\to Y$ can be
factorized as $j=ki$ as in
\[ \xymatrix{
    0 \ar[r] \ar[d] & W \ar[d] \\
    X \ar[r]^-{i} & X+W \ar[r]^-k & Y } \]
where $k$ is bijective on objects, $W$ is the discrete category on the
set of all objects of $Y$ not in the image of $j$,  and the square is a pushout. Since
the class $\cj$ contains $0\to W$,  is stable under pushouts, and
closed under composition, it will suffice to show that $\cj$ contains
the bijective on objects functors.

Suppose then that $j\colon X\to Y$ is bijective on objects, and
$p\colon A\to B$ is in $\cp$. In the diagram
\[ \xymatrix{
    A^Y \ar@/^1pc/[drr]^{A^j} \ar@/_1pc/[ddr]_{p^Y} \ar[dr]^{j\,\widetilde{\pitchfork}\,p}\\
    & P \ar[r]^-g \ar[d] & A^X \ar[d]^{p^X} \\
    & B^Y \ar[r]_-{B^j} & B^X } \]
where $P$ is a pullback, we are to show that
$j\,\widetilde{\pitchfork}\,p$ lies in $\cp$. But $B^j$ is a discrete
isofibration for any bijective on objects $j$, thus so is its pullback
$g$; likewise $A^j$ is a discrete isofibration. By the cancellation
property for discrete isofibrations, it follows that
$j\,\widetilde{\pitchfork}\,p$ is a discrete isofibration and so in $\cp$.
\end{proof}

\begin{example}
The representable fibrations, normal isofibrations, and discrete
isofibrations all satisfy the stability properties of
Lemma~\ref{lemma:Leibniz}.  Therefore each of these classes satisfies
the Leibniz condition with respect to the injective on objects functors.
\end{example}

\begin{remark}\label{remark:Groth}
  There are other interesting classes of representable isofibrations,
  such as the Grothendieck fibrations.  The lemma does not typically apply,
  since not every discrete isofibration is a Grothendieck fibration, and
  indeed the Leibniz condition does not typically hold. In the case $\ck=\Cat$,
for instance, consider the identity on objects inclusion $j\colon 2
\to \atwo$ from the discrete category with two objects to the generic
arrow $\{0 \to 1\}$,  and the Grothendieck fibration $p\colon \atwo
\to 1$. The Leibniz power of $p$ by $j$ is just $\atwo^j\colon\two^\atwo\to \atwo^2$, which can be seen as the inclusion of the
full subcategory of $\atwo^2=\atwo\x\atwo$ consisting of the three objects $(0,0)$,
$(0,1)$, and $(1,1)$. This is not a Grothendieck fibration: there is
no lifting of the map $(1,0)\to (1,1)$.
\end{remark}

\subsection{When all isofibrations are normal} \label{sect:NIP}

In this section we investigate the difference between representable
isofibrations and normal ones, as well as the case when this difference
disappears:

\begin{definition}
  A 2-category $\ck$ has the {\em normal isofibration property}, or
  {\em NIP}, if every representable isofibration in $\ck$ is a normal isofibration.
\end{definition}

We have already observed that $\Cat$ has this property.

Recall the construction of $w_f\colon A^\bbi\to L_f$ as in
Definition/Proposition~\ref{defprop:pseudolimits}.

\begin{proposition}\label{prop:f-vs-wf}
  For $f\colon A\to B$ 
  \begin{enumerate}
  \item $f$ is a representable isofibration if and only if $w_f$ is one
  \item $f$ is a normal isofibration if and only if $w_f$ is one.
  \end{enumerate}
In either case, $w_f$ is then a retract equivalence; if moreover $w_f$
is a normal isofibration, we call it a {\em normal retract
  equivalence}.
\end{proposition}

\begin{proof}
  Since $w_f$ is always an equivalence, it will be a retract
  equivalence if and only if it is a representable isofibration if and
  only if it has a section.
  
  To give a lifting $\kappa\colon x\cong u_f$ of $\lambda_f$ is
  equivalently to give a section of $w_f$. This proves (1). 

  A lifting $\kappa\colon x\cong u_f$ will have $\kappa d_f$ an
  identity if and only if and only if the corresponding section of
  $w_f$ is a filler for
  \[ \xymatrix{
      A \ar[r]^{\Delta} \ar[d]_{d_f} & A^\bbi \ar[d]^{w_f} \\
      L_f \ar[r]_-1 \ar@{.>}[ur] & L_f. } \]
  Since $d_f$ is an injective equivalence, there will be such a
  filler if $w_f$ is a normal isofibration. This gives one direction
  of (2).

  The converse follows by Lemma~\ref{lemma:Leibniz}, since $w_f$ is
  the Leibniz power of $f$ by $1\to\bbi$.
\end{proof}

\begin{corollary}
  The normal isofibration property holds just when the retract
  equivalences have the right lifting property with respect to the
  injective equivalences. 
\end{corollary}

\begin{proof}
Proposition~\ref{prop:f-vs-wf} implies that every representable isofibration is
normal just when this is the case for equivalences: in other words, when  every retract equivalence is a normal
isofibration; or, by Proposition~\ref{prop:wfs}, when the retract
equivalences have the right lifting property with respect to the
injective equivalences.
\end{proof}

We now specialize further to the case where $\ck$ is the 2-category
$\Cat(\ce)$ of internal categories in a complete category $\ce$. Then
$\ck$ is itself complete as a 2-category. Given an internal category
$A$ we write as usual $A_0$ for the object-of-objects, with a similar
notation for internal functors. We often identify an object of $\ce$
with the corresponding discrete internal category. 

\begin{proposition}
  The normal isofibration property holds in $\Cat(\ce)$ if and
  only if the split epimorphisms in $\ce$ have the right lifting
  property with respect to the split monomorphisms. 
\end{proposition}

\begin{proof}
  Consider a diagram as on the left 
  \[ \xymatrix{
      A \ar[r]^-c \ar[d]_i & C \ar[d]^p &
      A_0 \ar[r]^-{c_0} \ar[d]_{i_0} &  C_0  \ar[d]^{p_0}
      \\ B \ar[r]_-d & D &
    B_0 \ar[r]_-{d_0} \ar@{.>}[ur]^{h_0}  & D_0 } \]
  in which $i$ is an injective equivalence and $p$ a retract
  equivalence. Then we obtain a diagram in $\ce$ as on the right, with
  $i_0$ a split monomorphism and $p_0$ a split epimorphism. If this
  has a filler $h_0$, then since $i$ and $p$ are equivalences we may
  extend $h_0$ to a functor $h$ giving a filler in the diagram on the
  left.

  Conversely, suppose that we are given the solid part of the diagram
  on the right, with $p_0$ a split epimorphism and with $i_0$ a split
  monomorphism. If we make each object into a chaotic category with
  the given object-of-objects, we obtain a diagram as on the left,
  with $i$ an injective equivalence and with $p$ a retract
  equivalence. If this has a filler $h\colon B\to C$ then its action
  $h_0\colon B_0\to C_0$ on objects gives a filler in the diagram on
  the right. 
\end{proof}

\begin{corollary}
  If now $\ce$ is extensive, the normal isofibration property in
  $\Cat(\ce)$ is equivalent to the condition that every split monomorphism in
  $\ce$ is a coproduct injection.

  If $\ce$ is a topos, the normal isofibration property holds if and only if the topos $\ce$ is Boolean.
\end{corollary}

\begin{proof}
  The first part follows from
  \cite[Theorem~2.7]{RosickyTholen-FactorizationFibrationTorsion},
  which says that in an extensive category, the coproduct injections
  and split epimorphisms form a weak factorization system.

  A topos is Boolean when every monomorphism is a coproduct
  injection. But in a topos, every monomorphism is a pullback of a
  split monomorphism, and the coproduct injections are stable under
  pullback, so it suffices for every split monomorphism to be a
  coproduct injection.
\end{proof}

This provides a large class of examples which do {\em not} satisfy the
normal isofibration property. In particular,
$\ck=\Cat(\Set^\two)=\Cat^\two$ does not, since the topos $\Set^\two$ is not
Boolean.

Some examples of 2-categories, not of the type $\Cat(\ce)$, where the
normal isofibration property {\em does} hold can be found in
Examples~\ref{ex:TAlg} and~\ref{ex:VCat}, and in 
Section~\ref{sect:accessible}. 

\subsection{Preservation of isofibrations}

\begin{proposition}\label{prop:pres-isf}
  Suppose that the 2-category $\cl$ has pseudolimits of arrows, and
  the 2-functor $U\colon\cl\to\ck$ preserves them. Then $U$ preserves
  equivalences, injective equivalences, retract equivalences, normal
  retract equivalences,  representable isofibrations, and normal isofibrations.
\end{proposition}

\begin{proof}
  Any 2-functor at all preserves equivalences, injective equivalences,
  and retract equivalences. Preservation of representable
  isofibrations and normal isofibrations is an immediate consequence
  of Proposition~\ref{prop:characterizations-via-L}. Finally the
  normal retract equivalences are just the retract equivalences which
  are also normal isofibrations.
\end{proof}

\begin{proposition}\label{prop:refl-isf}
  Let $U\colon\cl\to\ck$ be as in Proposition~\ref{prop:pres-isf}, and
  suppose that each $U\colon \cl(A,B)\to\ck(UA,UB)$ is a discrete
  isofibration. Then $U$ also reflects representable
  isofibrations and normal isofibrations.

  Furthermore, the following are equivalent:
  \begin{enumerate}
  \item $U$ reflects equivalences
  \item $U$ reflects retract equivalences
  \item $U$ reflects normal retract equivalences
  \end{enumerate}
  and together they imply that
  \begin{enumerate}\addtocounter{enumi}{3}
  \item   $U$ reflects injective equivalences. 
  \end{enumerate}
\end{proposition}

\begin{proof}
  Let $f\colon A\to B$ in $\cl$ have pseudolimit $L_f$, with the usual
  notation. If $Uf$ is a representable isofibration, then there is a
  $\kappa\colon x\cong Uu_f$ with $Uf.\kappa=U\lambda_f$. We may lift this
  to a unique $\overline{\kappa}\colon y\cong u_f$ with
  $U\overline{\kappa}=\kappa$.

  Now $f\overline{\kappa}\colon fy\cong fu_f$ and $\lambda_f\colon v_f\cong
  fu_f$ satisfy
  \[ U(f\overline{\kappa})=Uf.U\overline{\kappa}=Uf.\kappa=U\lambda_f \]
  and so by the discrete isofibration property,
  $f\overline{\kappa}=\lambda_f$. Thus $f$ is a representable
  isofibration.

  If in fact $Uf$ is a normal isofibration, then we may choose
  $\kappa$ so that $\kappa.Ud_f=\id_{Uf}$. Then
  \[ U(\overline{\kappa}.d_f) = \kappa.Ud_f = \id_{Uf} = U\id_f \]
  and so $\overline{\kappa}.d_f=\id_f$, and $f$ is also a normal
  isofibration. This proves the first statement.

  Since the normal retract equivalences are the retract equivalences
  which are also normal isofibrations, and $U$ reflects normal
  isofibrations, (2) implies  (3).
  
  Suppose that $U$ reflects normal retract equivalences, and that $f\colon
  A\to B$ is such that $Uf$ is an equivalence. Then $Uv_f$ is a
  retract equivalence, thus $v_f$ is a normal retract equivalence, thus $f$
  is an equivalence. This proves that (3) implies (1).

  Now suppose that $U$ reflects equivalences, and that $f\colon A\to
  B$ is such that $Uf$ is an injective equivalence, say with
  retraction $r\colon UB\to UA$ and isomorphism $Uf.r\cong
  1$. Then $f$ is an (adjoint) equivalence, say with inverse
  equivalence $g\colon B\to A$, unit $\eta\colon 1\cong gf$, and
  counit $\epsilon\colon fg\cong 1$. Then $r.U\epsilon\colon
  Ug=r.Uf.Ug\cong r$ can be lifted to $\gamma\colon g\cong
  \overline{r}$, and now the isomorphisms $\gamma.f\colon g.f\cong
  \overline{r}.f$ and $\eta^{-1}\colon g.f\cong 1$ satisfy
  \[ U(\gamma.f) = U\gamma.Uf =
    r.U\epsilon.Uf=r.Uf.U\eta^{-1}=U\eta^{-1} \]
  and so $\gamma.f=\eta^{-1}$ and in particular
  $\overline{r}.f=1$, and $f$ is an injective equivalence. This
  proves that (1) implies (4). A similar argument shows that (1)
  implies (2).
\end{proof}

\begin{corollary}\label{cor:NIP}
  Suppose that the 2-category $\cl$ has pseudolimits of arrows, and
  the 2-functor $U\colon\cl\to\ck$ preserves them; and further that
  $U$ induces discrete isofibrations on the hom-categories. If $\ck$
  satisfies the normal isofibration property then so too does $\cl$. 
\end{corollary}

\begin{example}\label{ex:TAlg}
For a $2$-monad $T$ on $\ck$, let $T\textnormal{-Alg}$ be
the 2-category of strict $T$-algebras and pseudomorphisms, and
$U\colon T\textnormal{-Alg}\to\ck$ be the forgetful $2$-functor.  If $\ck$ has pseudolimits of arrows, then $U\colon T\textnormal{-Alg}\to\ck$ satisfies the assumptions of Corollary~\ref{cor:NIP}. Thus $T\textnormal{-Alg}$ satisfies the NIP if $\ck$ does so.
\end{example}

\begin{example}\label{ex:VCat}
  Let $\cv$ be a monoidal category with finite limits, and $\cv\Enr$ the 2-category of $\cv$-categories,
  $\cv$-functors, and $\cv$-natural transformations. Then the 2-functor
  $U\colon\cv\Enr\to\Cat$ sending a $\cv$-category to its underlying ordinary
  category satisfies the assumptions of Corollary~\ref{cor:NIP}; thus
  $\cv\Enr$ satisfies the normal isofibration property.
\end{example}

If in Corollary~\ref{cor:NIP}, we add the assumption that $U\colon\cl\to\ck$ is a fibration for the model structure on 2-$\Cat$ (of \cite{qm2cat}, but as corrected in \cite{qmbicat}), then the condition that $\cl$ has and $U$ preserves pseudolimits of arrows is in fact redundant.  This is the content of the following result.  

\begin{proposition}\label{prop:equiv-fibration}
  Let $\ck$ be a 2-category with pseudolimits of arrows, and
  $U\colon\cl\to\ck$ a 2-functor which
  \begin{enumerate}
  \item has the equivalence lifting property: if $A\in\cl$ and
    $u\colon L\to UA$ is an equivalence, then there is an equivalence
    $\overline{u}\colon \overline{L}\to A$ with $U\overline{u}=u$
  \item acts as a discrete isofibration $\cl(A,B)\to\cl(UA,UB)$ on hom-categories.
  \end{enumerate}
  Then $\cl$ has pseudolimits of arrows and $U\colon\cl\to\ck$
  preserves them. It then follows that if $\ck$ satisfies the normal
  isofibration property then so too does $\cl$. 
\end{proposition}

\begin{proof}
  Suppose that $f\colon A\to B$ is a morphism in $\cl$, and form
  the pseudolimit $L$ in $\ck$ of $Uf$, as displayed below on the
  left.  Then $u\colon L\to UA$ is an
  equivalence, so we may lift it to an equivalence
  $\overline{u}\colon \overline{L}\to A$. Now $\lambda\colon v\to
  Uf.u=U(f.\overline{u})$ is invertible, so there is a unique lifting
  $\overline{\lambda}\colon \overline{v}\cong f.\overline{u}$, as
  below in the centre.
  \[ \xymatrix @R1pc {
       & L \ar[ddl]_{u} \ar[ddr]^{v}
      \ar@{}[dd]|{\textstyle\cong_{\lambda}} &&&
      \overline{L} \ar[ddl]_{\overline{u}} \ar[ddr]^{\overline{v}} \ar@{}[dd]|{\textstyle\cong_{\overline{\lambda}}} &&&
      C \ar[ddl]_{a} \ar[ddr]^{b} \ar@{}[dd]|{\textstyle\cong_{\beta}}
      \\ 
      \\ 
      UA \ar[rr]_-{Uf} && UB & A \ar[rr]_-{f} && B & A \ar[rr]_-f && B } \]

  Now $u$ and $d$ are part of an adjoint equivalence with counit
  $ud=1$ and unit 
  $\delta\colon 1\cong du$. As in \cite[Proposition~6]{qmbicat}, we
  may lift this to an adjoint equivalence involving $\overline{u}$ and
  some $\overline{d}$ and unit $\overline{\delta}\colon 1\cong
  \overline{d}\overline{u}$, and by the uniqueness aspect of the discrete
  isofibration property the counit is once
  again an identity. Likewise, 
  $\overline{\lambda}\,\overline{d}\colon \overline{v}\overline{d}\cong
  f$ is a lifting of the identity $\lambda d\colon vd\cong f$ and so
  is itself an identity; in particular $\overline{v}\overline{d}=f$.

  We
  shall show that $\overline{u}$, $\overline{v}$, and
  $\overline{\lambda}$ exhibit $\overline{L}$ as the pseudolimit of
  $f$. By Definition/Proposition~\ref{defprop:equivalences}(3), it
  will suffice to check the 1-dimensional aspect of the universal
  property. 
  
  Suppose then that $a\colon C\to A$, $b\colon C\to B$ and
  $\beta\colon b\cong fa$, as displayed above on the right.  Apply $U$ to these data and use the
  universal property of $L$ to obtain a unique map $c\colon UC\to L$
  with $u.c=Ua$, $v.c=Ub$, and $\lambda.c=U\beta$.

  Now $\delta c\colon c\cong d.u.c=d.Ua=U(\overline{d}.a)$ has a
  unique lifting $\gamma\colon \overline{c}\cong \overline{d}a$. We
  shall show that $\overline{c}$ provides the required factorization
  through $\overline{L}$. Now
  \begin{itemize}
  \item $\overline{u}\gamma\colon \overline{u}\,\overline{c}\cong
    \overline{u}\overline{d}a=a$ is a lifting of an identity, thus an
    identity
  \item $\overline{\lambda}\overline{c}\colon
    \overline{v}\,\overline{c}\cong f\overline{u}\,\overline{c}=fa$
    is a lifting of $U\beta$, thus equal to $\beta$
  \end{itemize}
  and so $\overline{u}\,\overline{c}=a$,
  $\overline{v}\,\overline{c}=b$, and $\overline{\lambda}\overline{c}=\beta$.
  
  Thus we have the existence of a factorization. For the uniqueness
  aspect, suppose that $c'\colon C\to \overline{L}$ satisfies
  $\overline{u}c'=\overline{u}\overline{c}$,
  $\overline{v}c'=\overline{v}\overline{c}$, and
  $\overline{\lambda}c'=\overline{\lambda}\overline{c}$. Since
  $\overline{u}$ is an equivalence, there is a unique isomorphism
  $\theta\colon c'\cong \overline{c}$ with $\overline{u}\theta$ equal
  to the identity. Then $U\theta\colon c\to c$ satisfies $u.U\theta=1$
  and so, since $u$ is an equivalence, $U\theta=1$. Thus $\theta\colon
  c'\cong \overline{c}$ is a lifting of an identity, hence an
  identity, and $c'=\overline{c}$.

  This proves that $\cl$ has the pseudolimits and $U$ preserves
  them. The final sentence then follows by Corollary~\ref{cor:NIP}.
\end{proof}

\section{The $2$-cosmos structure on a 2-category with flexible limits}
\label{sect:flex}

In this section we give one of our main results, Theorem~\ref{thm:flexible}, which shows that any $2$-category with flexible limits forms a $2$-cosmos in which the chosen isofibrations are the normal isofibrations.  Moreover such $2$-cosmoi are \emph{Cauchy-complete}, in a sense to be defined.

Let us begin by reminding the readers about flexible limits in $2$-categories.  These were introduced in \cite{BKPS} and admit many equivalent formulations.  Perhaps the most elementary description is that they are those limits constructible from products, inserters, equifiers, and splittings of idempotents \cite{BKPS}, whilst a satisfying abstract description is that they are the cofibrantly-weighted $2$-categorical limits \cite{hty2mnd}.   

The characterization below, proven by the first-named author as Proposition A.1 of \cite{Bourke2019Accessible}, shows that that flexible limits and normal isofibrations are closely related.  We include a proof since the construction of pullbacks of normal isofibrations from flexible limits is important to our concerns.

\begin{proposition}\label{prop:Bourke-flexible}
  For a 2-category \ck, the following are equivalent
  \begin{enumerate}
  \item $\ck$ has flexible limits
  \item $\ck$ has products, powers, splittings of idempotents, and pullbacks of normal
    isofibrations
  \item $\ck$ has products, powers, splittings of idempotents, and pullbacks of discrete isofibrations.
  \end{enumerate}
\end{proposition}
\proof
Since discrete isofibrations are normal isofibrations, certainly (2)
implies (3). In order to prove that (3) implies (1), we use the fact
that flexible limits are generated by products, inserters, equifiers,
and splittings of idempotents. Thus it will suffice to show that (3)
implies the existence of inserters and equifiers. The inserter of
$f,g\colon A\to B$, and the equifier of $\beta,\beta'\colon h\to
k\colon A\to B$ can be constructed as the pullbacks below.
\[ \xymatrix{
    \operatorname{Ins}(f,g) \ar[r] \ar[d] & B^\atwo \ar[d] &
    \operatorname{Eq}(\beta,\beta') \ar[r] \ar[d] & B^{\atwo}
    \ar[d] \\
    A \ar[r]^{\binom{f}{g}} & B^2 & A \ar[r]^{\binom{\beta}{\beta'}} & B^{\two_2} } \]
Here $2\to \atwo$ is as in Remark~\ref{remark:Groth}, whilst $\two_2\to\two$ is the map from the generic parallel pair of arrows $\{0 \rightrightarrows 1\}$ to $\two$ which identifies the parallel arrows.  Since both functors are bijective on objects, it follows that the vertical restriction maps on the right of both squares are discrete isofibrations.

It remains to show that (1) implies (2), and since products, powers,
and splittings of idempotents are all flexible, we need only prove the
existence of pullbacks along normal isofibrations. This is the key
contribution of \cite[Proposition~A.1]{Bourke2019Accessible}.

Let $f\colon A\to C$ be a normal isofibration, and $g\colon B\to C$
arbitrary, and form the isocomma object
\[ \xymatrix {
    P \ar[r]^-p \ar[d]_q^{~}="1" & A \ar[d]^f_{~}="2" \\
    B \ar[r]_-g & C
    \ar@{=>}"1";"2"^{\phi} } \]
and the lifting $\xi\colon x\cong p$ with $f\xi=\phi$ (and so also $fx=gq$).

By the universal property of the isocomma object, there is a unique
$e\colon P\to P$ with $pe=x$, $qe=q$, and $\phi e=\id_{fx}$. Then $\phi ee$
is also an identity, $qee=qe$, and $pee=xe$; while $\xi e\colon
xe\cong pe$ is the lifting of the identity $\phi e\colon gq=gq=fx$
which must be the identity by normality. It follows by the universal
property of the isocomma object that $ee=e$, so that $e$ is
idempotent. By the 2-dimensional aspect of the universal property,
there is a unique isomorphism $\epsilon \colon e\cong 1$ with
$p\epsilon=\xi$ and $q\epsilon$ the identity. 

We claim that a splitting of $e$ gives the desired pullback of $f$ and
$g$. To show this, we should show that $ey=y$ if and only if $\phi y$
is an identity. On the one hand, if $ey=y$ then $\phi y=\phi ey$ which
is an identity since $\phi e$ is one. On the other hand, if $\phi y$
is an identity then so is its lifting $\xi y\colon xy=pey\cong py$,
and so in particular $pey=py$. We also know that $qey=qy$, since
$qe=q$; while $\phi ey$ is the identity on $fpey=fpy$, as is $\phi
y$. Thus $ey=y$ by the universal property of the isocomma object.
Note also that then $\epsilon y$ is the identity.

If $e$ splits as $e=ir$ with $ri=1$, then a straightforward
calculation shows that $pi$ and $qi$ give the desired pullback.
The two-dimensional aspect of the universal property follows easily
from the fact that the splitting of $e$ will be equivalent to $P$.
\endproof

We now turn to the (new) concept of Cauchy-completeness in the context of $\infty$-cosmoi and $2$-cosmoi.

\begin{definition}
An $\infty$-cosmos $\ck$ is {\em Cauchy-complete} if $\ck$ has split idempotents and any retract of an isofibration is an isofibration.   A $2$-cosmos $\ck$ is Cauchy-complete if the $\infty$-cosmos $n\ck$ is Cauchy complete. 
\end{definition}

Since the Cauchy-completeness condition for an $\infty$-cosmos does not involve the simplicial enrichment, clearly we have

\begin{proposition}
A $2$-cosmos $\ck$ is Cauchy complete just when it has split idempotents and a retract of an isofibration is an isofibration. 
\end{proposition}

With this, we are in a position to prove:

\begin{theorem}\label{thm:flexible}
Let $\ck$ be a $2$-category.  Then $\ck$ has flexible limits if and only if the normal isofibrations equip $\ck$ with the structure of a Cauchy complete $2$-cosmos.
\end{theorem}
\begin{proof}
One direction follows immediately from Proposition~\ref{prop:Bourke-flexible}.

For the converse, suppose that $\ck$ admits flexible limits. We shall
verify the conditions of Proposition~\ref{prop:2cosmos}.  Since each normal isofibration is a representable isofibration, Condition (a) holds.  

Since $\ck$ has flexible limits, it has products, powers, splittings
of idempotents and pullbacks of normal isofibrations by
Proposition~\ref{prop:Bourke-flexible} again.  As far as the existence of limits in Condition (b) goes, it remains to construct limits of towers of normal isofibrations.  

Consider a tower of normal isofibrations in $\ck$, with objects $A_n$ for
$n<\omega$, and morphisms $f^m_n\colon A_m\to A_n$ for $n\le
m<\omega$. Form the pseudolimit $P$ with projections $p_n\colon P\to
A_n$ and structure isomorphisms $\pi^m_n\colon p_n\cong f^m_n p_m$.

\subsubsection*{Step 1.} Define inductively a cone over the tower with
vertex $P$ and morphisms $q_n\colon
P\to A_n$, and with isomorphisms $\alpha_n\colon q_n\cong
p_n$, compatible in the sense that
\[ \xymatrix{
    f^m_nq_m \ar@{=}[r] & q_n \ar[r]^-{\alpha_n} & p_n \ar[r]^-{\pi^m_n} & f^m_n p_m } \]
is equal to $f^m_n\alpha_m$ for all $n<m$.

Start by setting $q_0=p_0$ and $\alpha_0=1$.

For a successor ordinal $n+1$, we have $p_{n+1}\colon P\to A_{n+1}$,
and isomorphisms
\[ \xymatrix{
    q_n \ar[r]^-{\alpha_n} & p_n \ar[r]^-{\pi^{n+1}_n} &
    f^{n+1}_np_{n+1} } \]
and so may use the normal isofibration structure to lift to an
isomorphism $\alpha_{n+1}\colon q_{n+1}\to p_{n+1}$ with
$f^{n+1}_nq_{n+1}=q_n$ and 
$f^{n+1}_n\alpha_{n+1}=\pi^{n+1}_n \alpha_n$.
A straightforward induction on $m$ shows that
$f^m_n\alpha_m=\pi^m_n.\alpha_n$ for all $m>n$. 

\subsubsection*{Step 2.} Next we show that a morphism $x\colon X\to P$
has $\alpha_n x\colon q_n x\to p_nx$ an identity for all $n<\omega$ if and
only if $\pi^m_n x\colon p_nx\to f^m_np_mx$ is an identity for all
$n<m <\omega$.

Since the diagram
\[ \xymatrix{
    q_nx \ar[r]^-{\alpha_n x} \ar@{=}[d] & p_n x \ar[d]^{\pi^m_n x} \\
    f^m_n q_m x \ar[r]_-{f^m_n\alpha_m x} & f^m_n p_m x} \]
commutes, if each $\alpha_n x$ is an identity then so is each $\pi^m_n
x$.

Suppose conversely that each $\pi^m_nx$ is an identity. We show by
induction that each $\alpha_n x$ is an identity.

First of all $\alpha_0$ is itself an identity hence so is $\alpha_0x$.

For a successor ordinal $n+1$, we obtain $\alpha_{n+1}x\colon
q_{n+1}x\to p_{n+1}x$ as a lift of
\[ \xymatrix{ q_n x \ar[r]^-{\alpha_n x} & p_n x \ar[r]^-{\pi^{n+1}_n
      x} & f^{n+1}_n p_{n+1}x } \]
so it will suffice, using normality, to show that this is an identity. But
$\alpha_n x$ is an identity by inductive hypothesis, and
$\pi^{n+1}_nx$ by assumption. 

\subsubsection*{Step 3.} Next we define an idempotent $e\colon P\to  
P$; later, we shall see that splitting $e$ gives the desired limit.

By the universal property of the pseudolimit, there is a unique
induced $e\colon P\to P$ with $p_n e=q_n$ for all $n<\omega$, and
\[ \xymatrix{
    q_n \ar@{=}[r] & p_n e \ar[r]^-{\pi^m_n e} & f^m_np_me \ar@{=}[r]
    & f^m_n q_m } \]
equal to the identity for all $n<m<\omega$. By Step~2 it follows that
$\alpha_n e\colon q_ne\to p_ne$ is an identity for all $n<\omega$;
and now $p_nee=q_ne=p_ne$ for all $n$, while each $\pi^m_nee$ and
$\pi^m_ne$ are both the (same) identity, whence $e^2=e$ by the
universal property of the pseudolimit $P$.

\subsubsection*{Step 4:} Next we obtain a (strict) cone over the
original tower.

Let the idempotent $e\colon P\to P$ split as
\[ \xymatrix{
    P \ar[r]^-r & L \ar[r]^-i & P } \]
and observe that
\[ f^m_n p_m i r = f^m_n p_m e = f^m_n q_m = q_n = p_n e = p_n i r \]
for all $n<m$, and so, cancelling the split epimorphism $r$, we deduce
that $f^m_n p_m i = p_n i$. Thus the $p_n i\colon
L\to A_n$ define a strict cone.

\subsubsection*{Step 5:} This cone is a limit.
If the maps $x_n\colon X\to A_n$
define another cone, then there is a unique $x\colon X\to P$ with
$p_n x=x_n$ for all $n$ and $\pi^m_n x$ an identity for all $n<m$. 
By Step~2, it follows that $\alpha_n x\colon q_nx\to p_nx$ is an
identity for all $n<\omega$, and so in particular that $q_n
x=p_nx$. Now $p_n ex=q_nx=p_nx$ for all $n$, while $\pi^m_n ex$ and
$\pi^m_nx$ are both the (same) identity, and so $ex=x$ by the
universal property of the pseudolimit $P$. It follows
that $p_nirx=p_nex=p_nx=x_n$ for all $n<\omega$, and so $rx$ is a factorization.

For the uniqueness, suppose that $y$ is any map for which $p_niy=x_n$ for all
$n<\omega$. Since $\pi^m_n i=\pi^m_niri=\pi^m_nei$ is an identity for
all $n<m<\omega$, so too is $\pi^m_n iy$, and thus $iy=x$ by the
universal property of the pseudolimit. Then $y=riy=rx$ giving the
uniqueness.

We leave the 2-dimensional aspect of the universal property as an
exercise for the reader.

Finally, we turn to the properties in Condition (c).  Clearly, each
map $p\colon A \to 1$ is a normal isofibration: indeed, since the only
isomorphism $pu \cong v \in \ck(A,1)$ is the identity, defining its
lifting to be $1\colon u \cong u$ gives a normal cleavage.

Stability of normal isofibrations under products, pullbacks and limits of towers follows from the fact, established in Proposition~\ref{prop:wfs}, that they are the right class of a weak factorisation system on $\ck$.  
Since the normal isofibrations are closed under composition and contain the discrete isofibrations, they satisfy the Leibniz property by Lemma~\ref{lemma:Leibniz}.
\end{proof}

The following result shows that normal isofibrations are the minimal choice of isofibration in the Cauchy-complete setting.

\begin{theorem}\label{thm:minimal}
In a $2$-cosmos, each normal isofibration is a retract of a chosen isofibration.  In particular, in a Cauchy-complete $2$-cosmos, each normal isofibration is chosen.
\end{theorem}
\proof
Let $f\colon A\to B$ be arbitrary. We may form the
pseudolimit of $f$ as the pullback
\[ \xymatrix{
    L \ar[r] \ar[d]_{\binom{u}{v}} & B^\iso \ar[d] \\
    A\x B \ar[r]_{f\x1} & B\x B } \]
of the chosen isofibration $B^\iso\to B\x B$, and so $\binom{u}{v}$ is
a chosen isofibration; but so too are the projections of $A\x B$,
hence so too are $u$ and $v$.

If now $f$ is a normal isofibration, then the filler $w$ in
\[ \xymatrix{
    A \ar[d]_d \ar[r]^1 & A \ar[d]^f \\
    L \ar[r]_v \ar@{.>}[ur]^w & B } \]
exhibits $f$ as a retract of $v$.
\endproof

Earlier we studied 2-categories  in which {\em every} representable isofibration is
normal, and we shall return to this topic in
Section~\ref{sect:accessible}. In relation to such 2-categories we
have the following observation:

\begin{corollary}\label{cor:uniqueness}
  If $\ck$ is a 2-category with flexible limits in which every
  representable isofibration is normal, then there is a unique
  Cauchy-complete 2-cosmos structure on $\ck$.
\end{corollary}

\begin{proof}
  By definition every chosen isofibration is a representable
  isofibation. By the theorem, every normal isofibration is a chosen
  isofibration. Thus we have implications
  \[ \textnormal{normal} \implies \textnormal{chosen} \implies
    \textnormal{representable} \]
  between the classes of isofibrations.  If all representable isofibrations are normal, the result follows. 
\end{proof}
\begin{examples}\label{ex:flexible}
As explained in \cite{BKPS}, many important $2$-categories of
categorical structures and their pseudomorphisms admit flexible
limits.  For instance, the $2$-categories
\textbf{Lex},~\textbf{Reg},~\textbf{Ex},~\textbf{Coh}, and \textbf{Pretop},
of small finitely complete categories, regular
categories, Barr-exact categories, coherent categories, and pretopoi all admit flexible limits, as does the $2$-category $\textbf{MonCat}$ of monoidal categories and strong monoidal functors.  

Therefore, by Theorem~\ref{thm:flexible}, each of these gives an
example of a $2$-cosmos with isofibrations the normal isofibrations.
In all of the above examples, the forgetful 2-functor to \textbf{Cat}
satisfies the conditions of Proposition~\ref{prop:equiv-fibration},
and so the normal isofibrations and representable isofibrations
coincide, and are just the morphisms whose underlying functor is an isofibration of categories.

In Proposition 2.6 of \cite{Riehl2021On}, Riehl and Wattal show that
the 2-category $\textbf{Icon}$ of bicategories,  normal
pseudofunctors and icons is a $2$-cosmos with chosen
isofibrations the representable ones, equally the normal ones.  This
is another instance of the canonical $2$-cosmos structure on a
2-category with flexible limits.   As pointed out in the proof of Proposition 2.6 of \cite{Riehl2021On}, the NIP holds in $\textbf{Icon}$.

Similarly the $2$-category of bicategories, (not necessarily normal)
pseudofunctors, and icons is a $2$-cosmos.  Indeed, following Section
6.3 of \cite{Bourke2021Accessible} it has flexible
limits\begin{footnote}{Pie limits were discussed earlier in
    \cite{2nerves}.}\end{footnote}, and so by
Theorem~\ref{thm:flexible} is a $2$-cosmos with isofibrations the
normal ones --- once again, the NIP is satisfied. 
\end{examples}

\begin{example}
  We end this section with a non-example. Consider the 2-category
  \textbf{sMonCat} of strict monoidal categories, strong monoidal functors,
  and monoidal natural transformations.  By
  Proposition~\ref{prop:refl-isf}, both the isofibrations and the normal isofibrations are those
  strong monoidal functors whose underlying functors are (necessarily
  normal) isofibrations of categories. We shall show that
  \textbf{sMonCat} does not have pullbacks of these (normal)
  isofibrations, and so that choosing them as isofibrations does not make
  \textbf{sMonCat} into  a 2-cosmos.
  
  By~\cite[Section~6.2]{Bourke2019Accessible}, idempotents do not
  split in \textbf{sMonCat}.  To see this, let $C$ be a small monoidal category.  Its strictification
  $A \in \textbf{sMonCat}$ has as objects the words in $C$, while its
  morphisms are such that the mapping  $q \colon A \to C$ which
  evaluates words using the left bracketing is fully faithful.  Then
  $q \colon A \to C$ is a retract equivalence in \textbf{MonCat}, with
  section $r \colon C \to A$ sending each object to the corresponding
  word of length 1,
  and so can be made into an adjoint equivalence with unit $\alpha \colon 1 \cong r q$.    In particular, $e = rq$ is then an idempotent in \textbf{sMonCat} and the triangle equations yield $e \alpha = id = \alpha e$.  If $e$ splits in \textbf{sMonCat}, then it follows that $C$ is isomorphic to the splitting, a strict monoidal category.  But not every monoidal $C$ is isomorphic to a strict monoidal category: as explained in \emph{loc.cit.}, an example is given by the skeletal category of at-most-countable sets.

So suppose that $e \colon A \to A$ as above does not split, and let
$M$ be the chaotic category with two objects $0$ and $1$, made into a
strict monoidal category under addition mod 2. There is a unique
functor $\chi\colon A\to M$ sending $a\in A$ to $0$ if $ea=a$ and to
$1$ otherwise. It admits a unique structure of strong monoidal functor
and is, moreover, an isofibration: on the one hand, if $\chi a=1$ then
$a\cong ea$ and $\chi ea=0$; on the other, if $\chi a=0$ then $a\cong
a.i$, where $i$ is the monoidal unit,  and $\chi (a.i)=1$ since $a.i$ is a not a singleton.


Suppose that there were a pullback in \textbf{sMonCat} as below left
 \[ \xymatrix{
    B \ar[r]^-f \ar[d] & A \ar[d]^\chi \\ 1 \ar[r]_-0 & M }
    \hspace{2cm}
    \xymatrix{
    A \ar[r]^-e \ar[d] & A \ar[d]^\chi \\ 1 \ar[r]_-0 & M }
\]
where $0\colon 1\to M$ is the (strict) monoidal
functor picking out the unit object $0$.  Since $e$ is idempotent, the
square above right commutes, so that we obtain a unique morphism
$k \colon A \to B$ to the pullback such that $fk = e$.  We claim that $kf = 1_B$.   Since $0 \colon 1 \to M$ has terminal domain, it is monic, whence so is its pullback $f$.  Therefore it suffices to show that $fkf = f$ or, since $fk = e$, that $ef = f$.  Indeed $\alpha_{fb} \colon fb \cong efb$ has identity component for each $b \in B$ since $fb$, being sent by $\chi$ to $0$, is a word of length $1$, and so in the image of $e$.  Hence $f = ef$.  Therefore $B$ splits $e$, a contradiction.
\end{example}

\section{Accessible $2$-cosmoi}\label{sect:accessible}

In this section, we investigate the notion of accessible $2$-cosmoi,
the $2$-categorical version of the accessible $\infty$-cosmoi
introduced in \cite{Bourke2021Accessible}.  In
Theorem~\ref{thm:variousconditions}, we show that the condition
concerning bicolimits in an accessible $2$-cosmos can be replaced by a
simpler condition concerning colimits of pointwise equivalences being
equivalences.  Combining this with the existing literature, we are
able to give many examples of accessible $2$-cosmoi including all of
the examples of Cauchy complete $\infty$-cosmoi in
Example~\ref{ex:flexible}.

To get started, let us recall the definition of accessible $\infty$-cosmos, before breaking down what it means.

\begin{definition}\label{def:acc}
An $\infty$-cosmos $\ck$ is said to be \emph{accessible} if 
\begin{enumerate}
\item $\ck$ is accessible as a simplicially-enriched category;
\item the chosen isofibrations are accessible;
\item there exists a regular cardinal $\lambda$ such that $\lambda$-filtered colimits exist in $\ck$ and are homotopy colimits.
\end{enumerate}
\end{definition}

First, by Proposition 2.4 of \cite{Bourke2021Accessible}, $\ck$ is accessible as a simplicially-enriched category just when $\ck_0$ is accessible and, for each simplicial set $X$, the powering functor $(X\pitchfork-)_0\colon\ck_0\to\ck_0$ is accessible.  

Since $\ck_0$ is accessible, so is the category of arrows $\ck^\two_0$.  We call a class of morphisms $\cj$ in $\ck$ \emph{accessible} just when the corresponding full subcategory of $\ck^\two_0$ is accessible and accessibly embedded in $\ck^\two_0$.  In particular, in Condition (2) we require that the class of chosen isofibrations be accessible.

With regards Condition (3), the colimit of a diagram $D\colon \cj \to \ck_0$ is said to be a homotopy colimit if the induced morphism $$p^*\colon (\cj^{op},\SSet)(\Delta 1,\ck(D-,X)) \to (\cj^{op},\SSet)(Q,\ck(D-,X))$$ is an equivalence of quasicategories, where $p\colon Q \to \Delta 1 \in [\cj^{op},\SSet]$ is a cofibrant replacement in the projective Joyal model structure.\begin{footnote}{Note that the property of being a homotopy colimit is independent of the choice of cofibrant replacement.}\end{footnote}

\begin{definition}
A $2$-cosmos $\ck$ is accessible if the $\infty$-cosmos $n\ck$ is accessible.
\end{definition}

The following proposition refers to bicolimits rather than homotopy colimits.  Bicolimits can be defined in the same way as homotopy colimits above, but with $p\colon Q \to \Delta 1 \in [\cj^{op},\Cat]$ a cofibrant replacement now in the projective model structure on $[\cj^{op},\Cat]$.

\begin{proposition}\label{prop:acc}
A $2$-cosmos $\ck$ is accessible if and only if
\begin{enumerate}
\item $\ck$ is accessible as a $\Cat$-enriched category;
\item the chosen isofibrations are accessible;
\item there exists a regular cardinal $\lambda$ such that $\lambda$-filtered colimits exist in $\ck$ and are bicolimits.
\end{enumerate}
\end{proposition}
\begin{proof}
As remarked above, $n\ck$ is accessible as a simplicially-enriched
category just when $\ck_0$ is accessible and, for each simplicial set
$X$, the powering functor $(X\pitchfork-)_0\colon\ck_0\to\ck_0$ is
accessible.  Since we have a natural isomorphism  $(X \pitchfork
A)_{n\ck} \cong (\Pi X \pitchfork A)_{\ck}$ and each small category is
of the form $\Pi X$ for some $X$, this is equally to say that the
powering functor $(X\pitchfork-)_0\colon\ck_0\to\ck_0$ is accessible
for each small category $X$.  But combined with accessibility of
$\ck_0$, this is exactly to say that $\ck$ is accessible as a
$\Cat$-enriched category.  Thus Condition (1) in Definition~\ref{def:acc}
corresponds to (1) above. There is no difference between the two Condition
(2)s: they involve only the underlying ordinary category.

It remains to prove the equivalence of Condition (3) with
Definition~\ref{def:acc}.(3).  For this, we will prove that a
(conical) colimit $\colim D$ in $\ck$ is a bicolimit just when it is a homotopy colimit in $n\ck$.

To do this, first observe that the Quillen adjunction $\Pi \dashv N$ gives rise to a Quillen adjunction $\Pi \circ - \dashv N \circ -$ between the projective model structures on $[\cj^{op},\Cat$] and on $[\cj^{op},\SSet]$.  Given a cofibrant replacement $p\colon Q \to \Delta 1 \in [\cj^{op},\SSet]$, we will show that $\Pi p\colon  \Pi Q \to \Pi 1 = \Delta 1$ is a cofibrant replacement in $[\cj^{op},\Cat]$.  Indeed $\Pi Q$ is cofibrant since $\Pi \circ - $ is left Quillen.  Moreover, at $X \in \cj$, $p_X$ is a weak equivalence of cofibrant objects (all simplicial sets being cofibrant); hence the left Quillen $\Pi$ sends it to a weak equivalence $\Pi(p_X)$ so that, in particular, $\Pi(p)$ is a weak equivalence in $[\cj^{op},\Cat]$, as required. 

Therefore
$\colim D$ is a bicolimit just when $$\Pi p^*\colon (\cj^{op},\Cat)(\Pi \Delta 1,\ck(D-,X)) \to (\cj^{op},\Cat )(\Pi Q,\ck(D-,X))$$ is an equivalence of categories or, equivalently, just when the induced map $N(\Pi p^*)$ on nerves is an equivalence of quasicategories.  
But it follows easily from Lemma~\ref{lemma:powers} that this is isomorphic to
\[p^*\colon (\cj^{op},\SSet)(\Delta 1,n\ck(D-,X)) \to (\cj^{op},\SSet )(\Pi Q,n\ck(D-,X))\]
and hence an equivalence of quasicategories just when $\colim D$ is a homotopy colimit in $n\ck$.
\end{proof}

Condition (3) about bicolimits in the above result can be replaced by an apparently simpler condition, as in the following result.  

\begin{proposition}\label{prop:variant}
A $2$-cosmos $\ck$ is accessible if and only if it is satisfies Conditions (1) and (2) of Proposition~\ref{prop:acc} above and
\begin{enumerate}
\item[(3)] there exists a regular cardinal $\lambda$ such that $\lambda$-filtered colimits of pointwise equivalences in $\ck$ are equivalences.
\end{enumerate}
\end{proposition}
\begin{proof}
First note that a $2$-cosmos satisfying Conditions (1) and (2) above is Cauchy-complete, and hence has flexible limits by Theorem~\ref{thm:flexible}.  In order to prove the claim, it will therefore suffice to show that given an accessible $2$-category with flexible limits, Condition (3) above is equivalent to $\lambda$-filtered colimits in $\ck$ being bicolimits.  We defer the proof of this non-trivial result to the appendix, where it appears in more general form as Proposition~\ref{prop:bicolimits}.
\end{proof}

Building on Proposition~\ref{prop:variant}, we can give several different formulations of accessible $2$-cosmoi, each of which avoids any mention of bicolimits.  This result is closely connected to \cite[Proposition~8.1]{Bourke2019Accessible}, although here the isofibrations are not necessarily the representable ones.

\begin{theorem}\label{thm:variousconditions}
Let $\ck$ be a $2$-cosmos for which $\ck$ is accessible as a $\Cat$-enriched category.  The following are equivalent.
\begin{enumerate}
\item $\ck$ is an accessible $2$-cosmos;
\item The chosen isofibrations and equivalences are accessible;
\item $\ck$ is a Cauchy complete $2$-cosmos and the trivial fibrations are accessible.
\end{enumerate}
\end{theorem}
\begin{proof}
By Proposition~\ref{prop:variant}, we have $(2 \implies 1)$.  For $(1 \implies 3)$, we must show that the trivial fibrations are accessible.  This is established in Theorem 6.3 of \cite{Bourke2021Accessible}.  Assuming $(3)$, the proof that the equivalences are accessible is exactly as in the proof of Proposition 6.4 of \emph{loc.cit.}  To prove that the isofibrations are accessible, we prove that there is a
pullback
\[ \xymatrix{
    \ck^\isf_0 \ar[r] \ar[d] & \TF(\ck)_0 \ar@{>>}[d] \\
    \ck^\two_0 \ar[r]_S & \ck^\two_0 } \]
where $S$ is accessible.  The result then follows from the
 Makkai-Par\'{e} limit theorem \cite{MakkaiPare} for accessible
categories --- see Proposition 2.3 of \cite{Bourke2021Accessible} for
a description of the relevant special case of pullbacks.

Here $S\colon\ck^\two_0\to\ck^\two_0$ is the functor sending $f$ to
$w_f$, as in Definition/Proposition~\ref{defprop:pseudolimits}. This
is accessible, since it is constructed using finite limits.

It remains to prove that $f$ is an isofibration if and only if
$w_f\colon A^\bbi\to L$ is a trivial fibration; but $w_f$ is always an equivalence, so this
says that $f$ is an isofibration if and only if $w_f$ is one. Since
$w_f$ is the Leibniz power of $f$ by the injective-on-objects $1\colon
1\to \bbi$, one direction is clear. Suppose conversely that $w_f$ is
an isofibration. The composite
of $p_f\colon L\to B^\bbi$ and $\dom\colon B^\bbi\to B$ is the isofibration appearing in
the Brown factorization of $f$ \cite[Lemma~1.2.19]{elements}. Thus the
composite $\dom.p_f.w_f=\dom.f^\bbi=f.\dom\colon A^\bbi\to B$ is also an
isofibration; but $f$ is a retract of this, hence it too is an
isofibration. 
\end{proof}

In order to give examples of accessible $2$-cosmoi, we will use the following result.

\begin{corollary}\label{cor:normal}
Let $\ck$ be a $2$-category with flexible limits in which all representable isofibrations are normal.   Then the normal $2$-cosmos structure on $\ck$ of Theorem~\ref{thm:flexible} is accessible just when $\ck$ is an accessible $2$-category and the retract equivalences in $\ck$ are accessible.  
\end{corollary}
\begin{proof}
Recall that the normal $2$-cosmos structure is Cauchy-complete and has
as chosen isofibrations the normal ones.  By our additional assumption, these are just the representable isofibrations.  Hence a morphism is a trivial fibration just when it is an equivalence and a representable isofibration, which amounts to being a retract equivalence.  The result then follows directly from Theorem~\ref{thm:variousconditions}.
\end{proof}

\begin{examples}\label{ex:acc2cosmoi}
In \cite{Bourke2019Accessible}, the first-named author described a
class $\LP_{\cm}$ of $2$-categories which are closely related to the
above.  A $2$-category $\ck$ belongs to $\LP_\cm$ just when it is
accessible, has filtered colimits and flexible limits, finite flexible
limits commute with filtered colimits in $\ck$, and retract
equivalences in $\ck$ are accessible. Now if $\ck$ belongs to
$\LP_\cm$ and all representable isofibrations are normal then, by
Corollary~\ref{cor:normal}, it follows that $\ck$ is an accessible 2-cosmos.

As shown in \cite{Bourke2019Accessible}, many $2$-categories of
categorical structures and their pseudo-morphisms belong to
$\LP_{\cm}$.  For instance, each of the following from
Examples~\ref{ex:flexible} belongs to $\LP_\cm$:
\begin{itemize}
\item \textbf{Lex} by \cite[Section~6.4]{Bourke2019Accessible}
\item \textbf{Reg},~\textbf{Ex} by
  \cite[Section~6.5]{Bourke2019Accessible}
\item \textbf{Coh},~\textbf{Pretop} by adapting the techniques of \cite[Section~6.5]{Bourke2019Accessible}
\item \textbf{MonCat}, \textbf{Icon} by
  \cite[Section~7]{Bourke2019Accessible}.
\end{itemize}
As observed in Examples~\ref{ex:flexible}, in each case the normal
isofibrations are just the representable isofibrations. It follows
that each of these examples is an accessible 2-cosmos.

For $T$ a $2$-monad on a $2$-category $\ck$, let $w \in \{s,p,l,c\}$ stand for ``strict/pseudo/lax/colax", and let  $\textnormal{w-}T\textnormal{-Alg}$ denote the $2$-category of $w$-T-algebras and pseudomorphisms, and $U \colon \textnormal{w-}T\textnormal{-Alg} \to \ck$ the forgetful 2-functor.\\
If $T$ is a filtered colimit preserving $2$-monad on $\ck \in \LP_\cm$, and $w \neq s$, then by Proposition 7.1 of \cite{Bourke2019Accessible} $\textnormal{w-}T\textnormal{-Alg}$ also belongs to $\LP_\cm$.  
The corresponding result where $w=s$ holds only under further assumptions, such as when $T$ is flexible --- see Theorem 7.2 of \cite{Bourke2019Accessible} and the discussion below it.  In all of these cases, Corollary~\ref{cor:NIP} applies to $U$, so that if $\ck$ satisfies the NIP so does $\textnormal{w-}T\textnormal{-Alg}$ which thereby forms an accessible $2$-cosmos.\\
For $\ca$ a small $2$-category, we can consider the $2$-monad $T$ on $\Cat^{ob\ca}$ with the simple formula $TX(a) = \Sigma_{b \in \ca}\ca(b,a) \times Xb$ --- see Section 6.6 of \cite{BKP} for a complete description.  Then the above setting specialises to capture each of the $2$-categories of lax/pseudo/oplax functors $\ca \to \Cat$ and pseudonatural transformations between them, so that each of these form accessible $2$-cosmoi.\\
As proven in Theorem 5.8 of \cite{Bourke2019Accessible}, the $2$-category of $2$-functors and pseudonatural transformations $\ca \to \Cat$ belongs to $\LP_\cm$ if $\ca$ is cellular, and again it follows that it is an accessible $2$-cosmos.  A key special case is when $\ca$ is the cellular $2$-category $\two$ --- then the corresponding accessible $2$-cosmos is the 2-category $\Ps(\two,\Cat)$ of arrows and pseudo-commutative squares.
\end{examples}

On the other hand, we have the following non-example, making clear
that the accessibility of $\ck$ alone is not enough:

\begin{proposition}
  There is no accessible $\infty$-cosmos structure on $\Cat^\two$. 
\end{proposition}

\begin{proof}
If there were, then by Theorem~\ref{thm:minimal} every normal
isofibration would be a chosen isofibration, and so in particular a
representable isofibration. By Theorem~\ref{thm:variousconditions},
there would be a regular cardinal $\alpha$ for which the trivial
fibrations are closed under $\alpha$-filtered colimits. Since every
normal retract equivalence is a trivial fibration, and every trivial
fibration is a retract equivalence, this would imply that every
$\alpha$-filtered colimit of normal retract equivalences is a retract
equivalence. We shall show that this is impossible.

Let $\alpha$ be an arbitrary regular cardinal.
For any set $Y$, let $S_Y$ be the set of all subsets of $Y$ of
cardinality strictly less than $\alpha$, seen as a chaotic
category, so that the functor $S_Y\to 1$ is a retract equivalence.
Let $P_Y$ be the category whose objects are pairs $(U\in S_Y,
x\in U)$, made into a category in such a way that the projection
$\pi_Y\colon P_Y\to
Y$ is a retract equivalence. Write $f_Y$ for the resulting morphism
\[ \xymatrix{
    P_Y \ar[r] \ar[d] & Y \ar[d] \\ S_Y \ar[r] & 1 } \]
in $\Cat^\two$.

Since $Y$ and $1$ are discrete, there are no non-trivial liftings, and
$f_Y$ is a normal isofibration. By construction the maps $P_Y\to Y$
and $S_Y\to 1$ are retract equivalences; the map $f_Y$ will itself be
a retract equivalence if and only if it has a section.

If $Y$ has cardinality less than $\alpha$, then $f_Y$ does have a
section, with $1\to S_Y$ picking out $Y\in S_Y$ and with $Y\to P_Y$
sending $y$ to $(Y,y)$. On the other hand if $Y$ has cardinality
$\alpha$, then $f_Y$ does not have a section: there is no way to
choose $U\subseteq Y$ with $|U|<\alpha$ containing every element of
 $Y$. Finally the fact that $f_Y$ is an $\alpha$-filtered colimit of
the $f_X$ with $|X|<\alpha$ gives the desired contradiction.
\end{proof}

Perhaps the way to think about this is that if one wants to model
arrows in $\Cat$ in an $\infty$-cosmological way, the ``correct''
thing to do is not use $\Cat^\two$ itself but rather the full
subcategory $\Cat^\isf$ of isofibrations. This is accessible, by
\cite[Proposition~4.1]{Bourke2021Accessible} and the fact that $\Cat$
is itself an accessible $\infty$-cosmos. Alternatively, in the $2$-categorical setting one can use the accessible $2$-cosmos $\Ps(\atwo,\Cat)$ of arrows and pseudo-commutative squares, as discussed in Examples~\ref{ex:acc2cosmoi}.\black

\appendix
\section{}

Our goal is to prove the following proposition.
\begin{proposition}\label{prop:bicolimits}
Let $\ck$ be an accessible $2$-category with flexible limits, and suppose further that it admits $\cj$-colimits for some small category $\cj$.  Then $\cj$-colimits are bicolimits in $\ck$ if and only if $\cj$-colimits of pointwise equivalences in $\ck$ are equivalences.
\end{proposition}
We will break the proof of the non-trivial direction into a few steps.
We fix a small category $\cj$ and an accessible 2-category $\ck$ with
flexible limits and $\cj$-colimits. We further fix a trivial fibration
$t\colon K\to \Delta 1$ with cofibrant domain in the projective model
structure on $[\cj\op,\Cat]$. 

For any functor
$S\colon\cj\to\ck$ and object $A\in\ck$, there is an induced functor 
\[    t^*\colon  [\cj\op,\Cat](\Delta1,\ck(S,A)) \to
  [\cj\op,\Cat](K,\ck(S,A)) \]
and we are to show that this is an equivalence.

\begin{proposition}\label{prop:R}
  There is an enriched
  monad $R$ on $[\cj,\ck]$ sending an object $T\colon\cj\to\ck$ to
  $RT\colon\cj\to\ck$, where $RTJ=KJ\pitchfork TJ$. The unit $r\colon
  T\to RT$ is induced by $t\colon K\to\Delta 1$. 
\end{proposition}

\begin{proof}
The 2-category $[\cj\op,\Cat]$ is monoidal with respect to
the product, and there is an enriched action
\[ \Phi\colon [\cj\op,\Cat]\op\x [\cj,\ck]\to [\cj,\ck] \]
given by
\[ \Phi(F,T)J = FJ \pitchfork TJ. \]
Every object of $[\cj\op,\Cat]$ has a unique comonoid structure; in
particular, $K$ does so, and thus
$\Phi(K,-)\colon[\cj,\ck]\to[\cj,\ck]$ has the structure of a monad.
\end{proof}

\begin{proposition}\label{prop:wla}
For each $S\colon\cj\to\ck$ there
exists a $QS\colon \cj\to\ck$ and $\eta\colon S\to RQS$ such that the
induced map
\[  \xymatrix{
    [\cj,\ck](QS,T) \ar[r]^-{R} & [\cj,\ck](RQS,RT) \ar[r]^-{\eta^*} &
    [\cj,\ck](S,RT) } \]
is a retract equivalence. 
\end{proposition}

\begin{proof}
If $\ck$ is accessible and has flexible limits, then the same is true
of $[\cj,\ck]$, and moreover $R\colon[\cj,\ck]\to[\cj,\ck]$ will be
accessible and preserve these limits. Thus it will have a weak left
adjoint by \cite[Theorem~8.10]{BourkeLackVokrinek}. This means that
that for each $S\colon\cj\to\ck$ there exist a $QS$ and an $\eta\colon
S\to RQS$ such
that the displayed map is a shrinking morphism in $\Cat$; that is, a
retract equivalence.
\end{proof}

\proof[Proof of Proposition~\ref{prop:bicolimits}]
First we suppose that $\cj$-colimits of pointwise equivalences are equivalences in $\ck$.  Given $S\colon\cj\to\ck$, 
consider the following diagram, which commutes up to isomorphism
\[ \xymatrix{
    [\cj\op,\Cat](\Delta 1,\ck(S,A)) \ar[r]^-{t^*} \ar[d] &
    [\cj\op,\Cat](K,\ck(S,A)) \ar[d] \\
    [\cj,\ck](S,\Delta A) \ar[r]^-{r_*} \ar@{=}[d]& [\cj,\ck](S,R\Delta A)
    \ar[d]^{\ref{prop:wla}} \\
    [\cj,\ck](S,\Delta A) \ar[r]^-{q^*} \ar[d] & [\cj,\ck](QS,\Delta A)
    \ar[d] \\
    \ck(\colim S,A) \ar[r]^-{q^*} & \ck(\colim QS,A) } \]
wherein the arrow marked \ref{prop:wla} is the equivalence guaranteed by
Proposition~\ref{prop:wla} and $q\colon QS\to S$ corresponds under
this equivalence to $r\colon S\to RS$, and where the remaining unnamed
arrows are the canonical isomorphisms. We are to show that the top
horizontal is an equivalence, but this will be the case if and only if
the bottom horizontal is one, and so if $q$ induces an equivalence
$\colim QS\to \colim S$ in $\ck$.

Since $\cj$-colimits of pointwise equivalences are equivalences in $\ck$, it will
suffice to show that $qJ\colon QSJ\to SJ$ is an equivalence for all
$J\in\cj$. This in turn says that for each $J\in\cj$ and $A\in\ck$ the top row in
the following diagram is an equivalence
\[ \xymatrix{
    \ck(SJ,A) \ar[r]^-{qJ^*} \ar[d] & \ck(QSJ,A) \ar[d] \\
    [\cj,\ck](S,\cj(-,J)\pitchfork A) \ar[r]^-{q^*}
    \ar[d]_{[\cj,\ck](S,\pi_2\pitchfork A)} &
    [\cj,\ck](QS,\cj(-,J)\pitchfork A) \ar[d]^{\ref{prop:wla}} \\
    [\cj,\ck](S,(K\x\cj(-,J))\pitchfork A) \ar[r] &
    [\cj,\ck](S,R(\cj(-,J)\pitchfork A)) 
  } \]
where once again the arrow marked \ref{prop:wla} is the equivalence
guaranteed by the Proposition of that number, and the unnamed arrows
are the canonical isomorphisms. Thus it will suffice to show that
$[\cj,\ck](S,\pi_2\pitchfork A)$ is an equivalence, which in turn will
be the case if $\pi_2\colon K\x\cj(-,J)\to\cj(-,J)$ is one.

Now this $\pi_2$ is a pullback of the trivial fibration $K\to\Delta1$,
hence is itself a trivial fibration; that is, a pointwise retract
equivalence. In general, an equivalence-inverse of a pointwise
equivalence in $[\cj\op,\Cat]$ might be only pseudonatural, however
here $\pi_2$ is a trivial fibration with cofibrant codomain, and so
has a section in $[\cj\op,\Cat]$, and is therefore an equivalence in
the 2-category $[\cj\op,\Cat]$. This completes the difficult direction of the proof.

On the other hand, if $\cj$-colimits are bicolimits in $\ck$, then their universal property with respect to pseudo-natural transformations implies that they they take pseudo-natural equivalences to equivalences in $\ck$.  But since each pointwise equivalence is a pseudo-natural equivalence, this implies the result.
\endproof


\begin{thebibliography}{10}
  

\bibitem{BKPS}
G.~J. Bird, G.~M. Kelly, A.~J. Power, and R.~H. Street.
\newblock Flexible limits for {$2$}-categories.
\newblock {\em J. Pure Appl. Algebra}, 61(1):1--27, 1989.

 \bibitem{BKP}
 R.~Blackwell, G.~M. Kelly, and A.~J. Power.
 \newblock Two-dimensional monad theory.
 \newblock {\em J. Pure Appl. Algebra}, 59(1):1--41, 1989.

\bibitem{Bourke2019Accessible}
John Bourke. 
\newblock  {Accessible aspects of 2-category theory},
 \newblock {\em Journal of Pure and Applied Algebra}, 225(3), 
 106519, March 2021. 
 
 \bibitem{Bourke2021Accessible}
John Bourke and Stephen Lack.
\newblock  {Accessible $\infty$-cosmoi},  {\em Journal of Pure and Applied Algebra}, 227(5), 
 107255, May 2023. 

\bibitem{BourkeLackVokrinek}
  John Bourke, Stephen Lack, and Luk\'as Vok\v r\'inek.
  \newblock {Adjoint functor theorems for homotopically enriched
    categories}, {\em Advances in Mathematics} 412, 108812, Jan 2023.
    
\bibitem{Gambino:hty2limits}
Nicola Gambino.
\newblock Homotopy limits for 2-categories.
\newblock {\em Math. Proc. Cambridge Philos. Soc.}, 145(1):43--63, 2008.

\bibitem{Garner:2DTypeTheory}
Richard Garner.
\newblock Two-dimensional models of type theory.
\newblock {\em Math. Structures Comput. Sci.}, 19(4):687--736, 2009.

\bibitem{qm2cat}
Stephen Lack.
\newblock A {Q}uillen model structure for 2-categories.
\newblock {\em $K$-Theory}, 26(2):171--205, 2002.

\bibitem{qmbicat}
Stephen Lack.
\newblock A {Q}uillen model structure for bicategories.
\newblock {\em $K$-Theory}, 33(3):185--197, 2004.

\bibitem{hty2mnd}
Stephen Lack.
\newblock Homotopy-theoretic aspects of 2-monads.
\newblock {\em J. Homotopy Relat. Struct.}, 2(2):229--260, 2007.

\bibitem{2nerves}
Stephen Lack and Simona Paoli.
\newblock $2$-nerves for bicategories.
\newblock {\em K. Theory}, 38(2):153-175, 2008.


\bibitem{MakkaiPare}
Michael Makkai and Robert Par{\'e}.
\newblock {\em Accessible categories: the foundations of categorical model
  theory}, volume 104 of {\em Contemporary Mathematics}.
\newblock American Mathematical Society, Providence, RI, 1989.


\bibitem{Riehl2015The-2-category}
E.~Riehl and D.~Verity. 
\newblock The 2-category theory of quasi-categories.
\newblock {\em Adv. Math.}, 280:549--642, 2015.

\bibitem{Riehl2021On}
E.~Riehl and M.~Wattal. 
\newblock On $\infty$-cosmoi of bicategories,
{\em La Matematica} 1: 740--764, 2022.

\bibitem{elements}
E.~Riehl and D.~Verity.
\newblock {\em Elements of $\infty$-category theory},  volume 194 of Cambridge Studies in Advanced Mathematics, Cambridge University Press, 2022.


\bibitem{RosickyTholen-FactorizationFibrationTorsion}
Ji\v{r}\'{\i} Rosick\'{y} and Walter Tholen.
\newblock Factorization, fibration and torsion.
\newblock {\em J. Homotopy Relat. Struct.}, 2(2):295--314, 2007.

\bibitem{Street1972Fibrations}
Ross Street.
\newblock Fibrations and Yoneda's lemma in a $2$-category, 
\newblock{In {\em Category Seminar, Sydney}, {\em Lecture Notes in Math. vol. 420},
Springer (1974) 104--133.}


\end{thebibliography}
\end{document}